\documentclass[english,10pt]{article}
\usepackage[T1]{fontenc}
\usepackage[latin9]{inputenc}
\usepackage{amsmath}
\usepackage{amssymb}

\usepackage{mathrsfs}
\date{}
\usepackage{babel}

\begin{document}

\title{Mordell-Weil growth for GL2-type abelian varieties \\
over Hilbert class fields of CM fields }

\author{David Hansen%
\thanks{david\_hansen@brown.edu%
}}

\date{June 10, 2010}
\maketitle
\begin{abstract}
Let $A$ be a modular abelian variety of $\mathrm{GL}_{2}$-type over
a totally real field $F$ of class number one. Under some mild assumptions,
we show that the Mordell-Weil rank of $A$ grows polynomially over
Hilbert class fields of CM extensions of $F$.
\end{abstract}

\section{Introduction and statement of results}

Let $E$ be an elliptic curve defined over a number field $F$. The
Mordell-Weil group $E(\overline{F})$ is one of the most mysterious
groups in arithmetic. By now, there are many theorems giving partial
or complete descriptions of $E(K)$ for various classes of extensions
$K/F$ which are abelian or nearly abelian. The growth of $E(K)$
as $K$ ascends through some sequence of nearly abelian extensions
tends to be controlled by root numbers:
\begin{itemize}
\item When $F=\mathbf{Q}$ and $K_{\infty}$ is the anticyclotomic $\mathbf{Z}_{p}$-extension
of a fixed imaginary quadratic field $K$, Vatsal and Cornut ({[}V{]},
{[}CV{]}) show that as one ascends up the cyclic layers of $K_{\infty}$,
the rank of $E$ is controlled entirely by the root number of $E/K$.
\item When $F=\mathbf{Q}$ and $K_{d}$ is the Hilbert class field of an
imaginary quadratic field $\mathbf{Q}(\sqrt{-d})$, Templier shows
that the rank of $E$ over $K_{d}$ is at least $\gg d^{\delta}$
for some small but positive fixed $\delta$, provided that $\varepsilon(E)=-\varepsilon(E\otimes\chi_{-d})$.
In fact, Templier has given two distinct proofs of this theorem: a
short proof {[}Te2{]} built on the Gross-Zagier theorem and equidistribution
theorems for Galois orbits, and an analytic proof {[}Te1{]} which
analyzes an average value of L-functions directly, using tools from
analytic number theory.
\end{itemize}
Our aim in this paper is to generalize Templier's analytic proof to
totally real base fields. This is not a triviality, and leads us to
solve an interesting auxiliary problem concerning the meromorphic
continuation and growth of a certain Dirichlet series. 

To state our theorems, we introduce a little notation. Let $F$ be
a totally real number field of degree $d$; we assume for simplicity
that $F$ has class number one. Let $f$ be a non-CM holomorphic Hilbert
modular form over $F$ whose weights are all even, with trivial central
character. For any $\alpha\in\mathscr{O}_{F}$ which is totally positive,
the field $F(\sqrt{-\alpha})$ is a CM extension. Write $H_{\alpha}$
for the Hilbert class field of this extension. Write $\chi_{\alpha}$
for the quadratic idele class character of $F$ cut out by $F(\sqrt{-\alpha})$.
Now, let $\mathscr{C}(\alpha)$ be the group of everywhere unramified
idele class characters of $F(\sqrt{-\alpha})$, or equivalently the
character group of the class group of $F(\sqrt{-\alpha})$. By the
Brauer-Siegel theorem, this group is of size $|\mathscr{C}(\alpha)|\gg_{\varepsilon}(\mathbf{N}\alpha)^{\frac{1}{2}-\varepsilon}$
as $\mathbf{N}\alpha\to\infty$. By quadratic automorphic induction,
any $\chi\in\mathscr{C}(\alpha)$ gives rise to a holomorphic Hilbert
modular form $\theta_{\chi}$ over $F$, of parallel weight one. The
root number of the Rankin-Selberg L-series $L(s,f\otimes\theta_{\chi})$
is $\pm1$, and is independent of $\chi$; it equals $\varepsilon(f)\varepsilon(f\otimes\chi_{\alpha})$.

Our main result is the following theorem.

\textbf{Theorem 1.1. }\emph{Notation as above, if $\alpha\in\mathscr{O}_{F}^{+}$
is such that $\varepsilon(f)\varepsilon(f\otimes\chi_{\alpha})=-1$,
then we have\[
\frac{1}{|\mathscr{C}(\alpha)|}\sum_{\chi\in\mathscr{C}(\alpha)}L'(\tfrac{1}{2},f\otimes\theta_{\chi})\gg L(1,\mathrm{sym}^{2}f)\log(\mathbf{N}\alpha)L(1,\chi_{\alpha})\]
as $\mathbf{N}\alpha\to\infty$, where the implied constant depends
only on $F$.}

In fact we give a precise asymptotic for this average; see Theorem
3.2 (when $F=\mathbf{Q}$, Theorem 3.2 is one of the main results
of {[}Te1{]}). The subconvexity results of {[}MV{]} give $L'(\tfrac{1}{2},f\otimes\theta_{\chi})\ll\mathbf{N}\alpha^{\frac{1}{2}-\delta}$
for some fixed positive $\delta$, so via the Brauer-Siegel theorem
we immediately deduce

\textbf{Corollary 1.2. }\emph{Notation and assumptions as above, there
exists some $\delta>0$ such that at least $\gg(\mathbf{N}\alpha)^{\delta-\varepsilon}$
of the central derivatives $L'(\tfrac{1}{2},f\otimes\theta_{\chi})$
are nonvanishing.}

The zeta function of $H_{\alpha}$ factors as\[
\zeta_{H_{\alpha}}(s)=\prod_{\chi\in\mathscr{C}(\alpha)}L(s,\theta_{\chi}).\]
Hence, feeding Corollary 1.2 into the results of {[}TZ{]} (see e.g.
Theorem 4.3.1 of {[}Zh{]}) yields

\textbf{Corollary 1.3. }\emph{Let $A/F$ be a modular abelian variety
of $\mathrm{GL}_{2}$-type, with associated Hilbert modular form $f$.
Then $\mathrm{rank}A(H_{\alpha})\gg\mathrm{dim}A\cdot(\mathbf{N}\alpha)^{\delta-\varepsilon}$
as $\mathbf{N}\alpha\to\infty$ along any sequence with $\varepsilon(f\otimes\chi_{\alpha})=-\varepsilon(f)$.}

We turn to an overview of the proof of Theorem 1.1. The first step
is to give an expression for $L'(s,f\otimes\theta_{\chi})$ as a short
Dirichlet polynomial essentially of length $\mathbf{N}\alpha$. Averaging
over $\chi$ yields an expression of shape\[
\frac{1}{|\mathscr{C}(\alpha)|}\sum_{\chi\in\mathscr{C}(\alpha)}L'(\tfrac{1}{2},f\otimes\theta_{\chi})\approx\sum_{x,y\in\left(\mathscr{O}_{F}\times\mathscr{O}_{F}\right)/\Delta U_{F}}\frac{\lambda_{f}(x^{2}+\alpha y^{2})}{(\mathbf{N}(x^{2}+\alpha y^{2}))^{\frac{1}{2}}}V\left(\mathbf{N}(x^{2}+\alpha y^{2})\right).\]
Here $\lambda_{f}$ are the Hecke eigenvalues of $f$, indexed by
ideals of $\mathscr{O}_{F}$, and $V(x)$ is a smooth function $\mathbf{R}_{>0}\to\mathbf{R}$
which decays rapidly for $x\gg(\mathbf{N}\alpha)^{1+\varepsilon}$
and diverges near the origin like $V(x)\sim c\log x$. Splitting off
the $y=0$ term yields a main term. Our problem then reduces to estimating
sums of the form\begin{equation}
\sum_{\gamma\in\mathscr{O}_{F}}\lambda_{f}(\gamma^{2}+\beta)W(\gamma^{2}+\beta)\end{equation}
for $\beta\in\mathscr{O}_{F}^{+}$ fixed and $W$ smooth, in various
ranges.

So far we have followed, in this reduction, Templier's analytic proof
{[}Te1{]} of Theorem 1.1 over $\mathbf{Q}$. Templier treats the sums
(1) over $F=\mathbf{Q}$ by a delicate and ingenious application of
the $\delta$-symbol method of Duke-Friedlander-Iwaniec, which is
in turn an elaboration of the original circle method of Hardy-Littlewood-Ramanujan.
Rather than trying to make the circle method work over number fields,
we analyze the sums (1) by the spectral theory of Hilbert modular
forms of half-integral weight. Our main result in this direction is
the following theorem.

\textbf{Theorem 1.4. }\emph{Notation as above, fix $\beta\in\mathscr{O}_{F}^{+}$
and define the Dirichlet series\[
\mathscr{D}_{f}(s;\beta)=\sum_{\gamma\in\mathscr{O}_{F}}\frac{\lambda_{f}(\gamma^{2}+\beta)}{(\mathbf{N}(\gamma^{2}+\beta))^{s}}.\]
Then $\mathscr{D}_{f}(s;\beta)$ admits a meromorphic continuation
to the whole complex plane. Furthermore, $\mathscr{D}_{f}(s)$ is
entire in the halfplane $\mathrm{Re}(s)>\frac{1}{4}$, with the exception
of at most finitely many poles in the interval $s\in[\tfrac{1}{4},\tfrac{1}{4}+\frac{\theta}{2}]$,
and satisfies the bound $\mathscr{D}_{f}(s;\beta)\ll e^{\pi d|s|}(1+|s|)^{A}\left(\mathbf{N}\beta\right)^{\frac{1}{2}-s-\frac{1}{16}(1-2\theta)}$
in that same halfplane.}

Actually we prove a slightly more general result dealing with general
quadratic polynomials, see Proposition 3.1. This result seems to be
new even over $\mathbf{Q}$. Here $\theta=\frac{7}{64}$ is the best
known bound towards the Ramanujan conjecture for $\mathrm{GL}_{2}$
over number fields {[}BB{]}. If we could only show this theorem with
the exponent $\frac{1}{2}-\frac{1}{16}(1-2\theta)$ replaced by $\frac{1}{2}$,
this would just barely fail to be strong enough to imply Theorem 1.1.
After giving a spectral expansion of $D_{f}(s;\beta)$, we eventually
deduce this crucial savings from a beautiful theorem of Baruch-Mao,
relating Fourier coefficients on $\widetilde{\mathrm{SL}_{2}}$ to
twisted L-values, which we control in turn via a subconvex bound due
to Blomer and Harcos {[}BH{]}.

This paper is organized as follows. In section two we review holomorphic
Hilbert modular forms and their L-functions, as well as the spectral
theory of Hilbert modular forms of half-integral weight. In section
three, we show that Theorem 1.1 is implied by Proposition 3.1, which
we prove in turn in section 4.

\subsubsection{Acknowledgements }

This material had its genesis in my 2010 Brown University senior honors
thesis {[}Ha{]}, where I proved Theorem 1.4 in the case $F=\mathbf{Q}$.
I would like to thank Jeff Hoffstein for many extremely helpful discussions
during the writing of my thesis, and for reading an earlier draft
of this paper. I am also grateful to Keith Conrad for several enlightening
conversations on the arithmetic of CM extensions, to Jordan Ellenberg
for answering a key question, and to Nicolas Templier for very helpful
comments on an earlier draft of this paper.

During the writing of this paper, I was supported by a Barry M. Goldwater
scholarship, and by the Josephine de K$\mathrm{\acute{a}}$rm$\mathrm{\acute{a}}$n
Fellowship Trust; it is a pleasure for me to acknowledge their generosity.

\section{Background and lemmas}

\subsection{Hilbert modular forms}

Fix $F/\mathbf{Q}$ totally real of degree $d$ and class number one.
Fix an ordering $\sigma_{1},\dots\sigma_{d}$ on the embeddings of
$F$ into $\mathbf{R}$. Write $\mathscr{O}_{F}$ for the ring of
integers of $F$, $\mathscr{O}_{F}^{+}$ for the totally positive
integers, $U_{F}$ for the unit group, and $U_{F}^{+}$ for the totally
positive units. We shall frequently use the fact that ideals in $\mathscr{O}_{F}$
are parametrized by the set $\mathscr{O}_{F}^{+}/U_{F}^{+}$. Fix
$\delta$ a totally positive generator of the different ideal of $F$.
Set $\Delta_{F}$ the absolute discriminant of $F$.

Let $g\in\mathrm{SL}_{2}(\mathscr{O}_{F})$ act on $z=(z_{1},\dots,z_{d})\in\mathfrak{H}^{d}$
in the usual way, i.e. via\[
\left(\begin{array}{cc}
a & b\\
c & d\end{array}\right)\cdot z=\left(\frac{\sigma_{1}(a)z_{1}+\sigma_{1}(b)}{\sigma_{1}(c)z_{1}+\sigma_{1}(d)},\dots,\frac{\sigma_{d}(a)z_{d}+\sigma_{d}(b)}{\sigma_{d}(c)z_{d}+\sigma_{d}(d)}\right).\]
For $\mathfrak{b},\mathfrak{c}\subset\mathscr{O}_{F}$, we define
congruence subgroups by\[
\Gamma_{0}(\mathfrak{b},\mathfrak{c})=\left\{ \left(\begin{array}{cc}
a & b\\
c & d\end{array}\right)\in\mathrm{SL}_{2}(\mathscr{O}_{F}),b\in\mathfrak{b},c\in\mathfrak{c}\right\} .\]
A Hilbert modular form of weight $(k_{1},\dots,k_{d})$ and level
$\mathfrak{a}$ is a function $f:\mathfrak{H}^{d}\to\mathbf{C}$ which
transforms as\[
f(\gamma z)=f(z)\cdot\prod_{j=1}^{d}(\sigma_{j}(c)z_{j}+\sigma_{j}(d))^{k_{j}}\]
for $\gamma=\left(\begin{array}{cc}
a & b\\
c & d\end{array}\right)\in\Gamma_{0}(1,\mathfrak{a})$. We shall restrict our attention to Hilbert modular forms of parallel
weight two. It will be abundantly clear at every step of the proof
that the general weight case is no harder, and in fact that we could
treat forms which are spherical at infinity if we so desired; we have
made this choice to simplify our notation.

A Hilbert modular form of parallel weight two has a Fourier expansion\[
f(z_{1},\dots z_{d})=\sum_{\alpha\in\mathscr{O}_{F}^{+}}\lambda_{f}(\alpha)(\mathbf{N}\alpha)^{\frac{1}{2}}e\left(\sigma_{1}(\delta^{-1}\alpha)z_{1}+\dots+\sigma_{d}(\delta^{-1}\alpha)z_{d}\right).\]
Here $\delta$ is a totally positive generator of the different ideal
of $F$, and the coefficients $\lambda_{f}(\alpha)$ depend solely
on the ideal generated by $\alpha$. The conductor of $f$ is the
unique ideal $\mathfrak{n}_{f}\subset\mathscr{O}_{F}$ such that $f$
is a new vector for $\Gamma_{0}(1,\mathfrak{n}_{f})$. Hereafter,
for any $\alpha\in F$ and any $z\in\mathfrak{H}^{d}$, we abbreviate
$\mathrm{tr}(\alpha z)=\sigma_{1}(\alpha)z_{1}+\dots+\sigma_{d}(\alpha)z_{d}$.
The L-function attached to $f$ has Dirichlet series\[
L(s,f)=\sum_{\alpha\in\mathscr{O}_{F}^{+}/U_{F}^{+}}\frac{\lambda_{f}(\alpha)}{\mathbf{N}\alpha^{s}}.\]

Given $\chi\in\mathscr{C}(\alpha)$, define \[
r_{\chi}(a)=\sum_{\mathfrak{a}\subset\mathscr{O}_{F(\sqrt{-\alpha})},\,\mathfrak{a}\overline{\mathfrak{a}}=(a)}\chi(\mathfrak{a}).\]
The function\[
\theta_{\chi}(z)=\sum_{a\in\mathscr{O}_{F}^{+}}r_{\chi}(a)e(\mathrm{tr}(\delta^{-1}az))\]
is a Hilbert modular form of parallel weight one and level $\Gamma_{0}(1,(\alpha))$
with central character $\chi_{\alpha}$. The Rankin-Selberg L-function
attached to $f\otimes\theta_{\chi}$ has Dirichlet series\[
L(s,f\otimes\theta_{\chi})=L^{(\mathfrak{n}_{f})}(2s,\chi_{\alpha})\sum_{a\in\mathscr{O}_{F}^{+}/U_{F}^{+}}\frac{\lambda_{f}(a)r_{\chi}(a)}{\mathbf{N}a^{s}}.\]
Here the superscript indicates removal of the Euler factors at primes
dividing $\mathfrak{n}_{f}$. We are assuming for simplicity that
$\alpha$ and $\mathfrak{n}_{f}$ are coprime. The completed L-function
is given by\[
\Lambda(s,f\otimes\theta_{\chi})=(\mathbf{N}\alpha\mathbf{N}\mathfrak{n}_{f})^{s}(2\pi)^{-2ds}\Gamma\left(s+\frac{1}{2}\right)^{d}\Gamma\left(s+\frac{3}{2}\right)^{d}L(s,f\otimes\theta_{\chi}).\]
This function satisfies $\Lambda(1-s,f\otimes\theta_{\chi})=\varepsilon(f)\varepsilon(f\otimes\chi_{\alpha})\Lambda(s,f\otimes\theta_{\chi}).$ 

We shall also require Hilbert modular forms of half-integral weight.
For a more detailed exposition of these see {[}Sh1{]} or {[}Ko{]}.
Set\[
\theta_{F}(z)=\sum_{\alpha\in\mathscr{O}_{F}}e\left(\mathrm{tr}(\delta^{-1}\alpha^{2}z)\right).\]
This is an anologue of Jacobi's theta function, and is a Hilbert modular
form of parallel weight $\frac{1}{2}$ for $\Gamma_{0}((2),(2))$;
c.f. {[}Sh1{]}. Define $j(\gamma,z)$ by \[
j(\gamma,z)=\frac{\theta_{F}(\gamma z)}{\theta_{F}(z)|cz+d|^{\frac{1}{2}}}.\]
 Note that $|j(\gamma,z)|=1$. Now, for $\nu=\frac{1}{2},\tfrac{3}{2}$
, and any ideal with $\mathfrak{n}\subset(2)$, define $\mathbf{H}_{\nu}(\mathfrak{n})$
to be the space of functions on $\mathfrak{H}^{d}$ which transform
as \[
\phi(\gamma z)=j(\gamma,z)^{2\nu}\phi(z),\,\forall\gamma\in\Gamma_{0}((2),\mathfrak{n}).\]
We could of course allow more general vector weights, but for our
purposes this is enough. This is a Hilbert space under the inner product\[
\left\langle \phi_{1},\phi_{2}\right\rangle =\int_{\Gamma_{0}(\mathfrak{n})\backslash\mathfrak{H}^{d}}\phi_{1}(z)\overline{\phi_{2}(z)}d\mu.\]
There is a large collection of commuting self-adjoint operators acting
on $\mathbf{H}_{\nu}(\mathfrak{n})$: the weight $\nu$ Laplace operator
$\Delta_{\nu}$ acts on each $z_{i}$-variable separately. Under the
action of these operators, $\mathbf{H}_{\nu}(\mathfrak{n})$ breaks
up as a direct sum of two orthogonal subspaces spanned by unitary
Eisenstein series and cusp forms, respectively. Suppose $\phi$ is
such that $\Delta_{\nu}^{(j)}\phi=\lambda_{\phi}^{(j)}\phi$ for $j=1,\dots,d$,
where the superscript indicates which $z_{i}$-variable we are acting
on. Define $t_{\phi}^{(j)}$ by $\lambda_{\phi}^{(j)}=\tfrac{1}{4}+(t_{\phi}^{(j)})^{2}$.
Then $\phi$ has a Fourier expansion \[
\phi(z)=\phi_{0}(y)+\sum_{\alpha\in\mathscr{O}_{F}\smallsetminus\{0\}}\frac{\rho_{\phi}(\alpha)}{\sqrt{|\mathbf{N}\alpha|}}\prod_{j=1}^{d}W_{\frac{\nu}{2}\mathrm{sign}(\sigma_{j}(\alpha)),it_{\phi}^{(j)}}(4\pi|\sigma_{j}(\delta^{-1}\alpha)|y_{j})e(\mathrm{tr}(\delta^{-1}\alpha x)).\]
If $\phi$ is a cusp form the term $\phi_{0}(y)$ vanishes.

Within each cuspidal $\left(\Delta^{(1)},\dots,\Delta^{(d)}\right)$-eigenspace
we take an orthonormal basis which furthermore is diagonalized for
all the $T_{\mathfrak{p}^{2}}$-Hecke operators, for all $\mathfrak{p}$
prime to $2\mathfrak{n}$. We then have the following crucial

\textbf{Lemma 2.1. }\emph{Under the above assumptions, the coefficients
$\rho_{\phi}(\alpha)$ satisfy the bound\[
\rho_{\phi}(\alpha)\ll\mathbf{N}\alpha^{\frac{1}{4}-\frac{1}{16}(1-2\theta)}\prod_{j=1}^{d}e^{\frac{\pi}{2}t_{\phi}^{(j)}}(1+|t_{\phi}^{(j)}|)^{A},\]
where $A$ is some large but fixed positive constant, and $\theta=\frac{7}{64}$
is the best known exponent toward the Ramanujan-Petersson conjecture
on $\mathrm{GL}_{2}$.}

\emph{Proof. }This follows from a theorem of Baruch-Mao {[}BM{]},
which gives a relation of the form\[
\frac{|\rho_{\phi}(\alpha)|^{2}}{\left\Vert \phi\right\Vert ^{2}}e(\phi)=\frac{L(\tfrac{1}{2},\Phi\otimes\chi_{\alpha})}{\left\Vert \Phi\right\Vert ^{2}}.\]
for $\alpha$ squarefree. Here $\Phi$ is the integral-weight Shimura
correspondant of $\phi$, and $e(\phi)$ is an archimedian integral,
which in our case can be computed as a ratio of $\Gamma$-functions.
The twisted L-values satisfy the bound $L(\tfrac{1}{2},\Phi\otimes\chi_{\alpha})\ll(\mathbf{N}\alpha)^{\frac{1}{2}-\frac{1}{8}(1-2\theta)}\prod_{j=1}^{d}(1+|t_{\phi}^{(j)}|)^{A}$
by the main theorem of {[}BH{]}. For $\alpha$ non-squarefree, the
$\rho_{\phi}(\alpha)$ can be expressed recursively in terms of their
values on squarefree divisors of $\alpha$.

We shall in fact require slightly more general theta functions than
$\theta_{F}$.

\textbf{Lemma 2.2. }\emph{For any $\beta\in\mathscr{O}_{F}$, the
function\[
\theta_{F}^{\beta}(z)=\sum_{\alpha\in\mathscr{O}_{F}}e\left(\mathrm{tr}\left(\delta^{-1}(\beta+2\alpha)^{2}z\right)\right)\]
is a Hilbert modular form of parallel weight $\tfrac{1}{2}$ and bounded
level.}

\[
\]

\section{Reduction to an analytic problem}

Let $f$ and $\theta_{\chi}$ be as in section 2. From now on we assume
that the global root number of $f\otimes\theta_{\chi}$ is $-1$.
Our goal is an asymptotic evaluation of the series\[
S=\frac{1}{|\mathscr{C}(\alpha)|}\sum_{\chi\in\mathscr{C}(\alpha)}L'(\tfrac{1}{2},f\otimes\theta_{\chi}).\]
Our first step is to derive an exact formula for $\Lambda'(\tfrac{1}{2},f\otimes\theta_{\chi})$.
Consider the integral\[
I=\frac{1}{2\pi i}\int_{(3)}\Lambda(s+\tfrac{1}{2},f\otimes\theta_{\chi})\cos\left(\frac{\pi s}{200}\right)^{-200}\frac{ds}{s^{2}}.\]
By construction the integrand has a simple pole at $s=0$ of residue
$\Lambda'(\tfrac{1}{2},f\otimes\theta_{\chi})$. On the other hand,
moving the contour of integration to $(-3)$ we derive\[
\Lambda'(\tfrac{1}{2},f\otimes\theta_{\chi})=(2\pi)^{-d}\sqrt{\mathbf{N}\alpha c_{f}}\cdot\sum_{(b,\mathfrak{n}_{f})=1}\chi_{\alpha}(b)\sum_{a}\frac{\lambda_{f}(a)r_{\chi}(a)}{\sqrt{\mathbf{N}a}\mathbf{N}b}V\left(\frac{\mathbf{N}ab^{2}}{\mathbf{N}\alpha\mathfrak{n}_{f}}\right),\]
where \[
V(x)=\frac{1}{4\pi i}\int_{(3)}\Gamma\left(s+1\right)^{d}\Gamma\left(s+2\right)^{d}\cos\left(\frac{\pi s}{200}\right)^{-200}(2\pi)^{-2ds}x^{-s}\frac{ds}{s^{2}}.\]
This gives\[
L'(\tfrac{1}{2},f\otimes\theta_{\chi})=\frac{1}{2}\sum_{b}\sum_{a}\frac{\lambda_{f}(a)r_{\chi}(a)}{\sqrt{\mathbf{N}a}\mathbf{N}b}V\left(\frac{\mathbf{N}ab^{2}}{\mathbf{N}\alpha c_{f}}\right).\]
Next we sum over $\chi\in\mathscr{C}(\alpha)$, giving by orthogonality\[
S=\frac{1}{\mathscr{C}(\alpha)}\sum_{\chi\in\mathscr{C}(\alpha)}L'(\tfrac{1}{2},f\otimes\theta_{\chi})=\frac{1}{2}\sum_{b\in\mathscr{O}_{F}^{+}/U_{F}^{+}}\sum_{\mathfrak{a}\subset\mathscr{O}_{F(\sqrt{-\alpha})},\,\mathfrak{a}\, principal}\frac{\lambda_{f}(\mathfrak{a}\overline{\mathfrak{a}})}{\sqrt{\mathbf{N}(\mathfrak{a}\overline{\mathfrak{a}})}\mathbf{N}b}V\left(\frac{\mathbf{N}(\mathfrak{a}\overline{\mathfrak{a}})\mathbf{N}b^{2}}{\mathbf{N}\alpha c_{f}}\right)\]
Let $\xi\in\mathscr{O}_{F(\sqrt{-\alpha})}$ be such that $\mathscr{O}_{F(\sqrt{-\alpha})}=\mathscr{O}_{F}+\mathscr{O}_{F}\xi$;
our class number assumption guarantees the existence of such a {}``relative
integral basis'' in our situation. Then the norms of principal ideals
in the ring $\mathscr{O}_{F(\sqrt{-\alpha})}$ are parametrized precisely
by the quadratic form $(\gamma+\xi\beta)(\gamma+\overline{\xi}\beta)=\gamma^{2}+(\xi+\overline{\xi})\beta+\xi\overline{\xi}\beta^{2}$,
where $\gamma$ and $\beta$ run over $\mathscr{O}_{F}\times\mathscr{O}_{F}$
modulo the diagonal action of the unit group. We now split the sum
$S$ into two sums, $S_{\mathrm{main}}$ and $S_{0}$ according to
whether $\beta=0$ or not. Thus we have\[
S=S_{\mathrm{main}}+S_{0},\]
with\[
S_{\mathrm{main}}=\frac{1}{2}\sum_{b\in\mathscr{O}_{F}^{+}/U_{F}^{+},\,(b,\mathfrak{n}_{f})=1}\chi_{\alpha}(b)\sum_{\gamma\in\mathscr{O}_{F}/U_{F}}\frac{\lambda_{f}(\gamma^{2})}{\mathbf{N}\gamma\mathbf{N}b}V\left(\frac{\mathbf{N}\gamma^{2}\mathbf{N}b^{2}}{\mathbf{N}(\alpha\mathfrak{n}_{f})}\right)\]
and

\[
S_{0}=\frac{1}{2}\sum_{b\in\mathscr{O}_{F}^{+}/U_{F}^{+},\,(b,\mathfrak{n}_{f})=1}\chi_{\alpha}(b)\sum_{\beta\in\left(\mathscr{O}_{F}\smallsetminus\{0\}\right)/U_{F}}\sum_{\gamma\in\mathscr{O}_{F}}\frac{\lambda_{f}(\gamma^{2}+(\xi+\overline{\xi})\beta+\xi\overline{\xi}\beta^{2})}{\sqrt{\mathbf{N}(\gamma^{2}+(\xi+\overline{\xi})\beta+\xi\overline{\xi}\beta^{2})}\mathbf{N}b}V\left(\frac{\mathbf{N}(\gamma^{2}+(\xi+\overline{\xi})\beta+\xi\overline{\xi}\beta^{2})\mathbf{N}b^{2}}{\mathbf{N}(\alpha\mathfrak{n}_{f})}\right).\]
To evaluate $S_{\mathrm{main}},$ note that it is given identically
by the contour integral\[
S_{\mathrm{main}}=\frac{1}{2\pi i}\int_{(3)}\frac{L^{(\mathfrak{n}_{f})}(2s+1,\chi_{\alpha})L(2s+1,\mathrm{sym}^{2}f)}{\zeta_{F}(4s+2)}\Gamma(s+1)^{d}\Gamma(s+2)^{d}\cos\left(\frac{\pi s}{200}\right)^{-200}(2\pi)^{-2ds}\left(\mathbf{N}\alpha\mathfrak{n}_{f}\right)^{s}\frac{ds}{s^{2}}.\]
Pushing the contour to $(-\tfrac{1}{4})$ we pick up a pole at $s=0$,
of residue \[
r(f,\alpha)=\frac{L^{(\mathfrak{n}_{f})}(1,\chi_{\alpha})L(1,\mathrm{sym}^{2}f)}{\zeta_{F}(2)}\left(\tfrac{1}{2}\log\mathbf{N}\mathfrak{n}_{f}+\tfrac{1}{2}\log\mathbf{N}\alpha+\frac{L'^{(\mathfrak{n}_{f})}}{L^{(\mathfrak{n}_{f})}}(1,\chi_{\alpha})+\frac{L'}{L}(1,\mathrm{sym}^{2}f)+c_{F}\right).\]
To estimate the integrand along the contour $\mathrm{Re}(s)=(-\tfrac{1}{4})$,
we use the subconvexity bound $L(\tfrac{1}{2}+it,\chi_{\alpha})\ll\left(\mathbf{N}\alpha\right)^{\frac{1}{4}-\frac{1}{16}(1-2\theta)}$.
Invoking the bound $|\zeta_{F}(1+it)|\gg_{F}\left(\log(|t|+3)\right)^{-d}$
({[}IK{]}, Ch. 5) we derive \[
S_{\mathrm{main}}=r(f,\alpha)+O\left(\mathbf{N}\alpha^{-\frac{1}{16}(1-2\theta)}\right).\]
Finally, we show the estimate $r(f,\alpha)\gg_{F,\varepsilon}L(1,\mathrm{sym}^{2}f)\left(\mathbf{N}\alpha\right)^{-\varepsilon}$.
To see this, we define a function on ideals in $\mathscr{O}_{F}$
by\[
\tau_{\alpha}(\mathfrak{a})=\sum_{\mathfrak{b}\supset\mathfrak{a}}\chi(\mathfrak{b}).\]
The generating function of this is simply\[
\sum_{\mathfrak{a}\subset\mathscr{O}_{F}}\tau_{\alpha}(\mathfrak{a})\mathbf{N}\mathfrak{a}^{-s}=\zeta_{F}(s)L(s,\chi_{\alpha}).\]
This immediately implies $\tau_{\alpha}(\mathfrak{a})\geq0$, and
if $\mathfrak{a}$ is a square (i.e. it has even valuation at every
finite place) then $\tau_{\alpha}(\mathfrak{a})=1$. Now we compute
the sum\[
S(X,\alpha)=\sum_{\mathfrak{a}\subset\mathscr{O}_{F}}\frac{\tau_{\alpha}(\mathfrak{a})}{\mathbf{N}\mathfrak{a}}\exp\left(-\frac{\mathbf{N}\mathfrak{a}}{X}\right)\]
in two different ways. On one hand, we may write\[
S(X,\alpha)=\frac{1}{2\pi i}\int_{(3)}\zeta_{F}(s)L(s,\chi_{\alpha})\Gamma(s-1)X^{s-1}ds;\]
 moving the contour to $\mathrm{Re}(s)=\tfrac{1}{2}$, we pass a pole
at $s=1$, and estimating the integral along $\mathrm{Re}(s)=\tfrac{1}{2}$
using the subconvex bound for $L(\tfrac{1}{2}+it,\chi_{\alpha})$
and the rapid decay of the gamma function yields the asymptotic \[
S(X,\alpha)=L(1,\chi_{\alpha})(\log X+\gamma_{F})+L'(1,\chi_{\alpha})+O(\mathbf{N}\alpha^{\frac{1}{4}-\frac{1}{16}(1-2\theta)}X^{-\frac{1}{2}}).\]
 On the other hand, $S(X,\alpha)\geq0$, so choosing $X=e^{-\gamma_{F}}\mathbf{N}\alpha^{\frac{1}{2}-\delta}$
with $0<\delta<\tfrac{1}{8}(1-2\theta)$ gives\[
(\tfrac{1}{2}-\delta)\log\mathbf{N}\alpha\cdot L(1,\chi_{\alpha})+L'(1,\chi_{\alpha})\geq0,\]
so \[
\tfrac{1}{2}\log\mathbf{N}\alpha+\frac{L'}{L}(1,\chi_{\alpha})\gg\log\mathbf{N}\alpha\]
and hence \begin{eqnarray*}
r(f,\alpha)\gg_{F}L(1,\mathrm{sym}^{2}f)L(1,\chi_{\alpha})\log\mathbf{N}\alpha.\\
\gg_{F,\varepsilon}L(1,\mathrm{sym}^{2}f)\left(\mathbf{N}\alpha\right)^{-\varepsilon}\end{eqnarray*}
by the Brauer-Siegel theorem. 

There are two notable cases where the clean estimate $r(f,\alpha)\gg_{F}L(1,\mathrm{sym}^{2}f)$
holds. If $F(\sqrt{-\alpha})$ does not contain a quadratic extension
of $\mathbf{Q}$, this follows immediately from Lemma 8 of {[}St{]}.
If $d>2$ and the Galois closure of $F$ has degree $d!$ over $\mathbf{Q}$,
this is a consequence of the beautiful results in {[}HJ{]}.

To bound the error term $S_{0}$, first complete the square, writing
\begin{eqnarray*}
\gamma^{2}+(\xi+\overline{\xi})\beta+\xi\overline{\xi}\beta^{2} & = & (\gamma+\tfrac{1}{2}(\xi+\overline{\xi})\beta)^{2}+(\xi\overline{\xi}-\tfrac{1}{4}(\xi+\overline{\xi})^{2})\beta^{2}\\
 & = & (\gamma+\tfrac{1}{2}(\xi+\overline{\xi})\beta)^{2}-\tfrac{1}{4}(\xi-\overline{\xi})^{2}\beta^{2}.\end{eqnarray*}
Recall that $\xi$ is a relative integral generator for the ring of
integers of a CM extension $F(\sqrt{-\alpha})/F$, and as such we
may write\[
\xi=\tfrac{1}{2}(x+y\sqrt{-\alpha})\]
for some $x,y\in\mathscr{O}_{F}$. Hence $\xi+\overline{\xi}=x$ is
an element of $\mathscr{O}_{F}$, and $\xi-\overline{\xi}=y\sqrt{-\alpha}$,
so $-\frac{1}{4}(\xi-\overline{\xi})^{2}=\tfrac{1}{4}y^{2}\alpha$
is a totally positive element of $F$, of absolute norm $\geq4^{-d}\mathbf{N}\alpha$.

Now, to evaluate $S_{0}$, we first execute the $b$-sum inside the
definition of $V$, giving\[
S_{0}=\sum_{\beta\in\left(\mathscr{O}_{F}\smallsetminus\{0\}\right)/U_{F}}\sum_{\gamma\in\mathscr{O}_{F}}\frac{\lambda_{f}((\gamma+\tfrac{1}{2}(\xi+\overline{\xi})\beta)^{2}-\tfrac{1}{4}(\xi-\overline{\xi})^{2}\beta^{2})}{\sqrt{\mathbf{N}((\gamma+\tfrac{1}{2}(\xi+\overline{\xi})\beta)^{2}-\tfrac{1}{4}(\xi-\overline{\xi})^{2}\beta^{2})}}W\left(\frac{\mathbf{N}((\gamma+\tfrac{1}{2}(\xi+\overline{\xi})\beta)^{2}-\tfrac{1}{4}(\xi-\overline{\xi})^{2}\beta^{2})}{\mathbf{N}(\alpha\mathfrak{n}_{f})}\right)\]
where we have set\begin{multline*}
W(x)=\sum_{b\in\mathscr{O}_{F}^{+}/U_{F}^{+},\,(b,\mathfrak{n}_{f})=1}\chi_{\alpha}(b)\left(\mathbf{N}b\right)^{-1}V(x\mathbf{N}b^{2})\\
=\frac{1}{4\pi i}\int_{(3)}L^{(\mathfrak{n}_{f})}(2s+1,\chi_{\alpha})\Gamma\left(s+\frac{1}{2}\right)^{d}\Gamma\left(s+\frac{3}{2}\right)^{d}\cos\left(\frac{\pi s}{200}\right)^{-200}(2\pi)^{-2ds}x^{-s}\frac{ds}{s^{2}}.\end{multline*}
The $W$-function satisfies the crude estimate $W(x)\ll_{A}x^{-A}$
for any fixed $A<100$.

The key ingredient in estimating $S_{0}$ is the following proposition,
which is a slight generalization of Theorem 1.4; we defer the proof
until section four.

\textbf{Proposition 3.1. }\emph{Let $P(x)=x^{2}+ax+b$ be a polynomial
with $a,b\in\mathscr{O}_{F}$ and $D=b-\tfrac{1}{4}a^{2}\in\tfrac{1}{4}\mathscr{O}_{F}$
totally positive. Then the Dirichlet series\[
\mathscr{D}_{f}(s;D)=\sum_{\gamma\in\mathscr{O}_{F}}\frac{\lambda_{f}(\gamma^{2}+a\gamma+b)}{\left(\mathbf{N}(\gamma^{2}+a\gamma+b)\right)^{s}}\]
admits a meromorphic continuation to the entire complex plane. Furthermore,
it is holomorphic in the half-plane $\mathrm{Re}(s)\geq\tfrac{1}{4}$
with the exception of at most finitely many poles in the interval
$[\tfrac{1}{4},\tfrac{1}{4}+\tfrac{\theta}{2}]$, and it satisfies
the bound\[
\mathscr{D}_{f}(s;D)\ll(1+|s|)^{A}e^{\pi d|s|}\left(\mathbf{N}D\right)^{\frac{1}{2}-s-\frac{1}{16}(1-2\theta)}\]
in that same half-plane, where the implied constant depends polynomially
on $f$.}

For $\beta$ fixed, the $\gamma$-sum is given exactly by the integral\[
I(\beta)=\tfrac{1}{2\pi i}\int_{(2)}\mathscr{D}_{f}(s+\tfrac{1}{2};\tfrac{1}{4}\alpha y^{2}\beta^{2})L^{(\mathfrak{n}_{f})}(2s+1,\chi_{\alpha})\Gamma\left(s+1\right)^{d}\Gamma\left(s+2\right)^{d}\cos\left(\frac{\pi s}{200}\right)^{-200}(2\pi)^{-2ds}\left(\mathbf{N}\alpha\mathfrak{n}_{f}\right)^{s}\frac{ds}{s^{2}}ds;\]
moving the contour to $\mathrm{Re}(s)=50$ is justified by the absence
of poles in that region, and the rapid decay of the integrand. To
estimate the integral along this contour we use the bound of Proposition
3.1, giving $ $ $\mathscr{D}_{f}(50+\tfrac{1}{2}+it,\tfrac{1}{4}\alpha y^{2}\beta^{2})\ll e^{\pi d|t|}\cdot\mathbf{N}(\alpha\beta^{2})^{-50-\frac{1}{16}(1-2\theta)}$.
The product of $\Gamma$-functions decays like $e^{-\pi d|t|}$ by
Stirling's formula, and the cosine decays like $e^{-\pi|t|}$, so
upon using the trivial bound $L(101+it,\chi_{\alpha})\asymp1$, the
integral converges absolutely and is bounded by $\left(\mathbf{N}\alpha\right)^{-\frac{1}{16}(1-2\theta)}\left(\mathbf{N}\beta\right)^{-100}$.
Inserting this bound into the definition of $S_{0}$, the $\beta$-sum
converges absolutely, giving \[
S_{0}\ll\mathbf{N}\alpha^{-\frac{1}{16}(1-2\theta)}.\]
Gathering results, we have proven

\textbf{Theorem 3.2. }\emph{Notation and assumptions as in Theorem
1.1, and assuming Proposition 3.1, we have\[
\frac{1}{|\mathscr{C}(\alpha)|}\sum_{\chi\in\mathscr{C}(\alpha)}L'(\tfrac{1}{2},f\otimes\theta_{\chi})=r(f,\alpha)+O_{F}(\left(\mathbf{N}\alpha\right)^{-\frac{1}{16}(1-2\theta)})\]
where\begin{eqnarray*}
r(f,\alpha) & = & \pi^{d}\frac{L^{(\mathfrak{n}_{f})}(1,\chi_{\alpha})L(1,\mathrm{sym}^{2}f)}{\zeta_{F}(2)}\left(\tfrac{1}{2}\log\mathbf{N}\mathfrak{n}_{f}+\tfrac{1}{2}\log\mathbf{N}\alpha+\frac{L'^{(\mathfrak{n}_{f})}}{L^{(\mathfrak{n}_{f})}}(1,\chi_{\alpha})+\frac{L'}{L}(1,\mathrm{sym}^{2}f)+c_{F}\right)\\
 & \gg_{F} & L(1,\mathrm{sym}^{2}f)L(1,\chi_{\alpha})\log\mathbf{N}\alpha.\end{eqnarray*}
}

\section{Sums of Hecke eigenvalues along quadratic sequences}

Fix $a,b\in\mathscr{O}_{F}$ with $D=b-\tfrac{1}{4}a^{2}$ totally
positive. Consider the integral\[
I_{f}(s;D)=\int_{\mathscr{O}_{F}\backslash\mathbf{R}^{d}\times(i\mathbf{R_{>0}})^{d}}|y|^{\frac{1}{4}}\theta_{F}^{a}(z)e(4D\delta^{-1}z)\overline{|y|f(4z)}|y|^{s}d\mu.\]
Inserting the Fourier expansions of $\theta_{F}^{a}$ and $f$ yields\begin{eqnarray*}
I_{f}(s;D) & = & \Delta_{F}\cdot\sum_{\gamma\in\mathscr{O}_{F}}\mathbf{N}(\gamma^{2}+a\gamma+b)^{\frac{1}{2}}\lambda_{f}(\gamma^{2}+a\gamma+b)\int_{(i\mathbf{R}_{>0})^{d}}e^{-16\pi\mathrm{tr}(\delta^{-1}(\gamma^{2}+a\gamma+b)y)}|y|^{s+\frac{1}{4}}d^{\times}|y|\\
 & =\Delta_{F}^{s+\frac{5}{4}} & (4\pi)^{-d(s+\frac{1}{4})}\Gamma(s+\tfrac{1}{4})^{d}\sum_{\gamma\in\mathscr{O}_{F}}\frac{\lambda_{f}(\gamma^{2}+a\gamma+b)}{(\mathbf{N}(\gamma^{2}+a\gamma+b))^{s-\frac{1}{4}}}.\end{eqnarray*}
On the other hand, the function $|y|^{\frac{5}{4}}f(4z)\overline{\theta_{F}^{b}(z)}$
is an automorphic form of weight $\tfrac{3}{2}$ for some arithmetic
group $\Gamma(b,\mathfrak{n}_{f})\subset\mathrm{SL}_{2}(\mathscr{O}_{F})$
whose index is bounded polynomially in $\mathbf{N}\mathfrak{n}_{f}$.
Write $\mathbf{H}_{\frac{3}{2}}(b,\mathfrak{n}_{f})$ for the group
of automorphic forms of weight $\frac{3}{2}$ and level $\Gamma(b,\mathfrak{n}_{f})$.
We may spectrally expand the function $|y|^{\frac{5}{4}}f(4z)\overline{\theta_{F}^{b}(z)}$
over an orthonormal basis of this space, giving \[
|y|^{\frac{5}{4}}f(4z)\overline{\theta_{F}^{b}(z)}=\sum_{\phi\in\mathbf{H}_{\frac{3}{2}}^{cusp}(b,\mathfrak{n}_{f})}\left\langle f,\theta_{F}^{b}\phi\right\rangle \phi(z)+\sum_{\gamma\in cusps}\int_{\mathbf{R}}\left\langle f,\theta_{F}E_{\gamma}(\bullet,\tfrac{1}{2}+it)\right\rangle E_{\gamma}(z,\tfrac{1}{2}+it)dt.\]
This spectral expansion converges absolutely, and thus we may insert
it into $I(s)$ and interchange the order of integration and summation,
giving\begin{eqnarray*}
I_{f}(s;D) & = & \Delta_{F}(\mathbf{N}D)^{-\frac{1}{2}}\sum_{\phi\in\mathbf{H}_{\frac{3}{2}}^{cusp}(4\mathfrak{n}_{f})}\frac{\left\langle f,\theta_{F}\phi\right\rangle }{\left\langle \phi,\phi\right\rangle }\rho_{\phi}(4D)\int_{(i\mathbf{R}_{>0})^{d}}|y|^{s-1}e^{-2\pi\mathrm{tr}(\delta^{-1}4Dy)}\prod_{j=1}^{d}W_{\frac{3}{4},it_{\phi}^{(j)}}(4\pi|\sigma_{j}(\delta^{-1}4D)|y_{j})d^{\times}|y|\\
 &  & +\mathrm{Eis.}\end{eqnarray*}
Strictly speaking, there is a contribution from weight-$\tfrac{3}{2}$
single-variable theta functions as well, but for $\theta$ such a
function, the inner product $\left\langle f,\theta_{F}^{b}\theta\right\rangle $
is a linear combination of nonzero multiples of $\mathrm{res}_{s=1}L(s,\mathrm{sym}^{2}f\otimes\eta)$
for $\eta$ some finite-order Hecke characters (cf. {[}Sh2{]}), and
any twist of $L(s,\mathrm{sym}^{2}f)$ is entire since we are assuming
that $f$ is not CM. By the Mellin transform formula\[
\int_{0}^{\infty}e^{-2\pi y}W_{\alpha,\beta}(4\pi y)y^{s}d^{\times}y=(4\pi)^{-s}\frac{\Gamma(s+\frac{1}{2}-\beta)\Gamma(s+\frac{1}{2}+\beta)}{\Gamma(s+1-\alpha)},\]
the $d$-fold integral evaluates to\[
\left(\mathbf{N}D\right)^{1-s}(4\pi)^{-ds}\Gamma(s-\tfrac{3}{4})^{-d}\prod_{j=1}^{d}\Gamma(s+it_{\phi}^{(j)}-\tfrac{1}{2})\Gamma(s-it_{\phi}^{(j)}-\tfrac{1}{2}).\]
Comparing these two expansions yields\[
\mathscr{D}_{f}(s;D)=C\cdot\Delta_{F}\cdot\left(\mathbf{N}D\right)^{\frac{1}{4}-s}\sum_{\phi\in\mathbf{H}_{\frac{3}{2}}^{cusp}(b,\mathfrak{n}_{f})}\frac{\rho_{\phi}(4D)\left\langle f,\theta_{F}\phi\right\rangle }{\left\langle \phi,\phi\right\rangle }\prod_{j=1}^{d}\frac{\Gamma(s+it_{\phi}^{(j)}-\frac{1}{4})\Gamma(s-it_{\phi}^{(j)}-\tfrac{1}{4})}{\Gamma(s-\frac{1}{2})\Gamma(s+\frac{1}{2})}+\mathrm{Eis}(s).\]
We shall only treat the cuspidal part of this expansion, the Eisenstein
terms being a great deal simpler. 

By Stirling's formula $\Gamma(\sigma+it)\asymp|t|^{\sigma-\frac{1}{2}}\exp(-\frac{\pi}{2}|t|)$,
the individual quotients of gamma functions for $s=\sigma+it$ are
bounded away from their poles by\[
(1+|t|)^{A}\exp\left(-\frac{\pi}{2}(|t+t_{\phi}^{(j)}|+|t-t_{\phi}^{(j)}|-2|t|)\right).\]
Using the identity $|a+b|+|a-b|-2|a|=2\max(|b|-|a|,0)$, the $d$-fold
product of gamma functions is bounded by \[
(1+|t|)^{A}\exp\left(-\pi\sum_{j=1}^{d}\max(|t_{j}|-|t|,0)\right).\]
The triple product $\left\langle f,\theta_{F}\phi\right\rangle $
is bounded as $\phi$ varies, by Cauchy-Schwarz, so using Lemma 2.1
we find\[
\mathscr{D}_{f}(s;D)\ll\left(\mathbf{N}D\right)^{\frac{1}{2}-s-\frac{1}{16}(1-2\theta)}\sum_{\phi\in\mathbf{H}_{\frac{3}{2}}^{cusp}(b,\mathfrak{n}_{f})}(1+|t|)^{A}\exp\left(\pi\sum_{j=1}^{d}\left(\frac{|t_{j}|}{2}-\max(|t_{j}|-|t|,0)\right)\right)\]
Only $\phi$'s with $|t_{\phi}^{(j)}|\leq2t$, $j=1..d$ contribute
to this sum, and the number of such eigenvalues is bounded polynomially
in $|t|$ and $\mathbf{N}\mathfrak{n}_{f}$ by a weak form of Weyl's
law, so summing their contribution trivially we conclude\[
\mathscr{D}_{f}(s;D)\ll\left(\mathbf{N}D\right)^{\frac{1}{2}-s-\frac{1}{16}(1-2\theta)}(1+|s|)^{A}e^{\pi d|s|},\]
away from the poles of the quotients of gamma functions, which occur
at the points $s=\frac{1}{4}\pm it_{\phi}^{(j)}$. By the Shimura
correspondence and {[}BB{]}, the numbers $t_{\phi}^{(j)}$ lie in
$\mathbf{R}\cup i[-\frac{\theta}{2},\frac{\theta}{2}]$. This concludes
the proof of Proposition 3.1.

The exponential factor $e^{\pi d|s|}$ appearing in the bound of Proposition
3.1 can likely be removed with a little more work. The key to doing
this would be to prove a triple product bound\[
\left\langle f,\theta_{F}\phi\right\rangle \ll\prod(1+|t_{\phi}^{(j)}|)^{A}e^{-\frac{\pi}{2}|t_{\phi}^{(j)}|}.\]
Here is one possible way to show this bound. Let $\pi$ and $\sigma$
be unitary automorphic representations of $\mathrm{PGL}_{2}/F$, not
both Eisenstein, and let $\widetilde{\sigma}$ be an automorphic representation
of $\widetilde{\mathrm{SL}_{2}}/F$ which lifts to $\sigma$ under
the Shimura-Shintani-Waldspurger correspondence. Let $\chi$ be a
quadratic idele class character, with associated one-variable theta
function $\theta_{\chi}$, and let $\pi_{\chi}$ be the automorphic
representation of $\widetilde{\mathrm{SL}_{2}}$ generated by the
adelic lift of $\theta_{\chi}$. Choose factorizable vectors $\varphi_{\pi}\in\pi,\,\varphi_{\widetilde{\sigma}}\in\widetilde{\sigma},\,\varphi_{\chi}\in\pi_{\chi}$,
and let $S$ denote the set of places where at least one of the three
local vectors is ramified. Then we conjecture a formula of the form\[
\frac{\left|\left\langle \varphi_{\widetilde{\sigma}}\varphi_{\chi},\varphi_{\pi}\right\rangle \right|^{2}}{\left\langle \varphi_{\widetilde{\sigma}},\varphi_{\widetilde{\sigma}}\right\rangle \left\langle \varphi_{\pi},\varphi_{\pi}\right\rangle }=C\frac{L(\tfrac{1}{2},\sigma\otimes\chi\otimes\mathrm{sym}^{2}\pi)}{L(1,\mathrm{ad}\pi)L(1,\mathrm{ad}\sigma)}\prod_{v\in S}\beta_{v}(\varphi_{\widetilde{\sigma},v},\varphi_{\pi,v},\varphi_{\chi,v}).\]
Here the L-functions appearing on the right are \emph{completed} with
their archimedian gamma factors, and the $\beta_{v}$'s are local
integrals. This is a simultaneous generalization of two important
formulas of Shimura: when $\pi$ arises from an Eisenstein series,
the period integral on the left was key in Shimura's original lifting
construction from half-integral weight forms to integral weight forms;
and when $\sigma$ arises from an Eisenstein series, this is the integral
representation for the symmetric square L-function given in {[}Sh2{]}.
For certain very special pairs $\pi,\sigma$ this conjecture is in
fact a theorem of Ichino {[}Ich{]}, and it seems quite reasonable
to adapt his technique for a general proof. Anyway, assuming this
formula, the convexity bound for L-functions combined with the exponential
decay of the archimedian gamma factors yields an immediate proof of
the purported triple product bound.

\section*{References}

{[}BM{]} E. Baruch and Z. Mao, \emph{Central values of automorphic
L-functions}, Geom. and Func. Anal. 17\\
\\
{[}BB{]} V. Blomer and F. Brumley, \emph{On the Ramanujan conjecture
over number fields}, preprint\\
\\
{[}BH{]} V. Blomer and G. Harcos, \emph{Twisted L-functions over
number fields and Hilbert's eleventh problem}, Geom. and Func. Anal.,
to appear\\
\\
{[}CV{]} C. Cornut and V. Vatsal, \emph{Nontriviality of Rankin-Selberg
L-functions and CM points}\\
\\
{[}CP-SS{]} J. Cogdell, I. Piatetski-Shapiro and P. Sarnak, \emph{Estimates
on the critical line for Hilbert modular L-function and applications},
unpublished\\
\\
{[}Ha{]} D. Hansen, \emph{Some new results in analytic number theory,
}2010 Brown University senior honors thesis\\
\\
{[}HJ{]} J. Hoffstein and N. Jochnowitz, \emph{On Artin's conjecture
and the class number of certain CM fields. I, II}, Duke Math. J. 59\\
\\
{[}Ich{]} A. Ichino, \emph{Pullbacks of Saito-Kurokawa lifts, }Invent.
Math. 162\\
\\
{[}IK{]} H. Iwaniec and E. Kowalski, \emph{Analytic Number Theory,
}AMS Colloquium Publications\\
\\
{[}Ko{]} H. Kojima, \emph{On the Fourier coefficients of Hilbert-Maass
wave forms of half integral weight over arbitrary algebraic number
fields}, J. Num. Theory 107\\
\\
{[}MV{]} P. Michel and A. Venkatesh, \emph{The subconvexity problem
for $\mathrm{GL}_{2}$}, Publ. Math. de l'IHES, to appear\\
\\
{[}Sh1{]} G. Shimura, \emph{On the Fourier coefficients of Hilbert
modular forms of half-integral weight}, Duke Math. J. 71\\
\\
{[}Sh2{]} G. Shimura, \emph{On the holomorphy of certain Dirichlet
series}, Proc. London Math. Soc. 31\\
\\
{[}St{]} H. Stark, \emph{Some effective cases of the Brauer-Siegel
theorem, }Invent. Math. 23\\
\\
{[}Te1{]} N. Templier, \emph{A non-split sum of Fourier coefficients
of modular forms}, Duke Math. J., to appear\\
\\
{[}Te2{]} N. Templier, \emph{Minoration du rang des courbes elliptiques
sur les corps de classes de Hilbert}, Compositio Math., to appear\\
\\
{[}TZ{]} Y. Tian and S. Zhang, \emph{Kolyvagin systems of CM points
on Shimura curves, }in preparation\\
\\
{[}V{]} V. Vatsal, \emph{Uniform distribution of Heegner points,
}Invent. Math. 150\\
\\
{[}Zh{]} S. Zhang, \emph{Barcelona notes on the arithmetic of Shimura
curves}
\end{document}